# A Train Formation Plan with Elastic Capacity for Large-Scale Rail Networks

Boliang Lin*

School of Traffic and Transportation, Beijing Jiaotong University, Beijing 100044, People's Republic of China

**Abstract**: This paper studies the problem of optimizing the train formation plan and traffic routing (TFP&TR) simultaneously. Based on the previous research of TFP&TR with determinate parameters, we consider the fluctuation of flows and associate the elastic features of yard capacities and link capacities. To demonstrate such elastic peculiarities, the fuzzy theory is introduced and the membership function is designed to designate the satisfaction degree for the volumes of traffic flows. A non-linear integer programming model is developed considering various operational requirements and a set of capacity constraints, including link capacity, yard reclassification capacity and the maximal number of blocks a yard can be formed, while trying to minimize the total costs of accumulation, reclassification, and transportation.

**Keywords**: Railway network, Traffic routing, Train formation plan, Fuzzy theory

## 1 Introduction

For railway freight flow organization, an essential problem is how to deliver shipments on a capacitated physical network optimally, that is, to determine the best traffic routing for each shipment and assign each shipment into the most reasonable freight train services. The state-of-the-art research in this paper is that optimizes the train formation plan (TFP) and traffic routing (TR) simultaneously. TFP determines which train services(blocks) should be provide between which pairs of yards and which shipment should be picked up by each train service. In China, TFP is viewed as a tactical problem, which has a two-year intermediate planning horizon. In the process of updating TFP, the shipments for a day are assumed to be given in advance, and the number of the shipments is determined according to the empirical data which is valued by highly-experienced service designers (Lin (2012)). The volume of shipments usually fluctuates from day to day, which means that the elastic capacities are required for yards and links. In addition, the development of yards and links will increase the corresponding capacities. As mentioned above, TFP is usually updated every two years. During this period, the facilities and equipment of yards and links may also develop and update, which makes the increase of transportation capacity of links and capacity (reclassification capacity and number of available tracks) of yards. This also shows that the elastic capacities constraints should be considered when updating TFP. Therefore, solving the TFP and TR(TFP&TR) problem with elastic feature of practical

---

* E-mail address: bllin@bjtu.edu.cn (B.-L. Lin).

significance and theoretically challenging.

In terms of traffic routing problem, Lansdowne (1981) considered the tracks that controlled by several carriers in a manner and developed an algorithm for routing freight over a rail network. The objective function of their work was to route the freight to minimize the number of interline transfers and to maximize the revenue division for the originating carrier. Crainic et al. (1984) studied the problem which includes routing freight traffic, scheduling train services, and allocating classification work among yards on a rail network. Haghani (1989) developed a nonlinear 0-1 programming model to describe the train routing, makeup, and empty car distribution model, combined with the characteristic of the model, a heuristic decomposition algorithm is designed to solve the model. For train formation plan, it is to determine which train service should be provided and which shipment should be picked up by each train service. As a result, the solution of the optimization for TFP can not only ensure the economical delivery of shipments and get better use of locomotives and railcars, but also coordinate the utilization of capacity of rail lines and yards, fully explore the potential capacity of railway transportation, improve transport efficiency, and reduce congestion (Mu et al. (2011), Lin (2017), Crainic (2000)). Lin (2012) presented a formulation and solution for the TFP, and a bi-level programming model is developed to described the problem. Combining the characteristic of the model, the paper deigned a simulated annealing algorithm and use a real-world example based on the China railway system to test the model and algorithm. Chen (2018) studied the TFP in the railway freight transportation industry, it is necessary to mention that the input to this problem in this paper is a railway network and a set of shipments, the paper developed a linear binary programming model to minimize the total sum of train accumulation costs and car classification costs, while considering the unitary rule and the intree rule based on the Chinese railway background, a novel solution methodology that applies a tree-based decomposition algorithm is designed to solve the model.

There have been researches on the integrated optimization of TFP&TR, for example, Lin (2021). This paper proposed a non-linear 0-1 programming model to address the integrated problem, and it aimed to minimize the total costs of railway operation. The capacity constraints of yards and links are rigid constraints in the model, while not considered the elastic feature. Based on this research, we improved the model and develop the elastic capacity constraints of yards and links.

## 2 Problem description

The problem that this paper studies is integrated optimization of TR and TFP, this section first describes the TR problem and then demonstrates the TFP problem using a railway train service network. Finally, the need for elastic capacity of yards and links is also analysed.

### 2.1 Traffic routing

We use a simple example to illustrate the traffic routing of shipments, which is shown in figure 1. Figure 1 contains six yards and six links. Some configurations regarding this network are as follows: the capacity of path B→C→D is 200 cars while the capacity of path B→F→D is 1000 cars, and that of A→B, D→E are also 1000 cars. Assuming that the volume of all shipments $N_{ij}$ is 100 cars except that the volume of $N_{BC}$ is 200 cars.

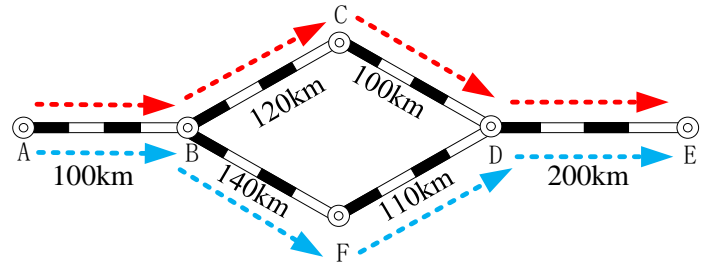

**Figure 1**: a simple railway network

Taking the shipment from A to E as an example, there are two path for $N_{AE}$ can be chosen: A→B→C→D→E(520km) and A→B→F→D→E(550km). In general, the shipment will be transported A→B→C→D→E which is the shortest path. However, the transportation capacity of link should be considered when studying the problem. The capacity of path B→C→D is 200 cars, and the volume of OD pair $N_{BC}$ is also 200 cars. This means that the transportation capacity of section B→C can only be used for $N_{BC}$, the other shipments whose origin is A and B only transported by B→F. Based on this, the transportation path for $N_{AE}$ is A→B→F→D→E. In the railway operation, this situation is very common and just only few shipments can be transported by the shortest path.

## 2.2 Train formation plan

The TFP is to determine which train service should be provided and which shipment should be picked up by each train service. In this section, we use the path A→B→C→D→E based on the above network to demonstrate the TFP.

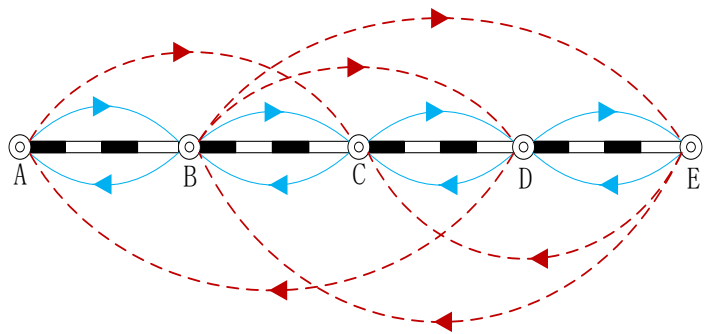

**Figure 2**: a simple-scale train service network

In figure 2, there are 14 train services in the network, including 6 direct train services represented by the red dashed line and 8 local train services for every adjacent yard represented by the blue line. It should be noted that the example is just used to illustrate the TFP, the direct train services in the network are generated randomly.

Besides, the capacity of yards is assumed to be sufficient. Based on the above train service network, the first phase of TFP is determined, that is, which train service should be provided. The second phase of TFP is to determine which shipment should be picked up by each train service. The result of this is to form a train service sequence for each shipment and determinee the shipping strategy for each shipment. For each shipment, it may have many shipping strategies, and the best result is that each shipment adopts the optimal shipping strategy to minimize the total operating cost, which is the motivation of our research.

Taking the shipment from A to E as an example, there are 4 shipping strategies for this shipment based on the above train service network, which are shown in figure 2. The possible shipping strategies for this shipment are as follows:

Strategy(a): A→B→C→D→E
Strategy(b): A→B→D→E
Strategy(c): A→B→E
Strategy(d): A→C→D→E

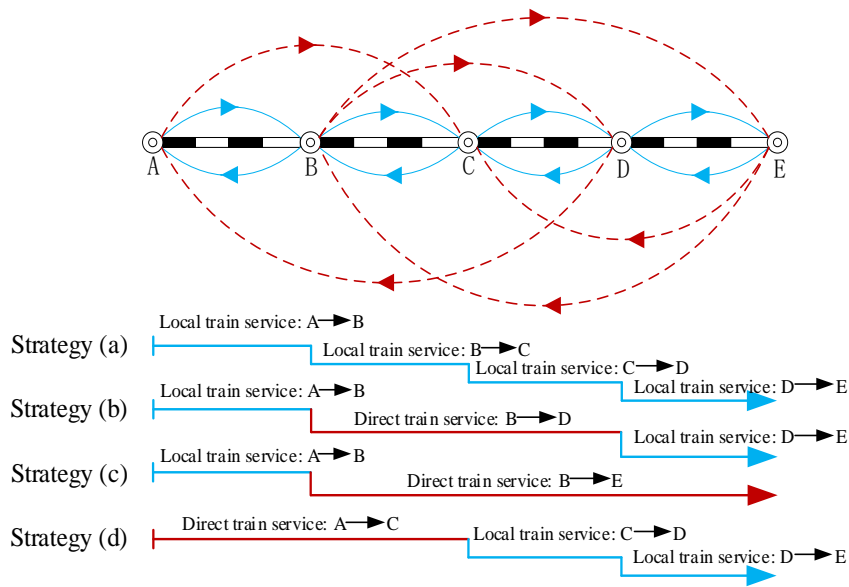

**Figure 3**: shipping strategies in the rail network

In strategy(a), the shipment is carried by the local train services without direct train services, and it needs to be reclassified at B, C,D. Different from strategy(a), the shipment is carried by the local train service A→B first, then merged with the shipment from B to E and carried by the direct train service B→D, then reach E by the local train service D→E in strategy(b). Compared with strategy(b), if the shipment is carried by strategy(c), it just needed to be transferred at B, and only has one reclassification operation. In strategy(d), the shipment is carried by the direct train service A→C first, after reclassified at yard C, it takes the local train service C→D and D→E to its destination. In summary, the second phase of TFP is to determine which shipping strategy for each shipment should be chosen by optimization analysis.

## 2.3 The elastic capacities of yards and links

We use the shipments A→E and C→E based on figure 2 as examples to illustrate the elastic reclassification capacity of yards. Besides, we assume that the above train service network is generated by such input data: the number of shipments from A to E is 100 cars per day, the number of shipments from C to E is 200 cars per day, the other shipments are 150 cars per day, the reclassification capacity of yard C is 300 cars per day, and the reclassification capacity of other yards is sufficient. If the optimal train service sequences for A→E and C→E are A→C→D→E, C→D→E respectively.

In the above example, the optimal train service sequence is generated according to the input data which includes the number of shipments per day. Here, we only analyse the shipments A→E and C→E, ignoring other shipments. The number of shipment A→E is 100 cars and it need to be reclassified at C according to the train service sequence. The number of shipment C→E is 200 cars, these two shipments will be integrated at yard C, the number of this merged shipment is 300 cars which is equal to the reclassification capacity of C. In the example, the number of daily shipments is regarded as a fixed value when performing optimization calculations, that is, assuming the number of shipments is the same every day. Actually, the number of shipments fluctuates every day, for example, the number of shipment A→E is 100 cars on the first day and 200 cars on the second day. If it is 200cars, the above train service sequence is obviously unreasonable, because the number of shipment at C is 400 cars(A→E:200 cars, C→E:200 cars) which beyond the reclassification capacity of C. Based on the above analysis, the storage capacity should be required for yards to adapt to fluctuations of shipments, that is, the elastic characteristic of reclassification capacity should be considered when studying TFP. In the same way, it is the same for the available tracks of yards and transportation capacity of links, we omit the narrative of it here. In addition, TFP is a tactical problem and usually updated every two years. During this period, the facilities and equipment of yards and links may also develop and updated, which makes the increase of them. This also shows that the elastic capacity should be considered when updating TFP.

## 3 Mathematical model

As mentioned in our previous work, we improve the traditional TFP&TR model which is referred in Lin (2021) and develop the elastic capacity constraints based on it. In this section, we first describe the traditional TFP&TR model, then introduce the improvements that include the elastic capacity constraints of yards and links. In addition, we use the membership function to discuss the satisfaction degree for the volumes of traffic flows.

In reference Lin (2021), the simultaneous optimization of TFP&TR problem can be formulated as a non-linear binary programming model whose objective function and constraints are expressed as follows:

$$\min \ Z(X) = \sum_{i \in V} \sum_{j \in V, j \neq i} c_i m_{ij} y_{ij} + \sum_{k \in P(i,j)} \tau_S^k \sum_{i \in V} \sum_{j \in V, j \neq i} f_{ij}(X) x_{ij}^k + \lambda \sum_{i \in V} \sum_{j \in V, j \neq i} D_{ij}(X) \sum_{l \in \rho(i,j) - \rho_{ij}^S} L_{ij}^l \xi_{ij}^l \quad (1)$$

*Subject to*

$$y_{ij} + \sum_{k \in P(i,j)} x_{ij}^k = 1, \quad \forall i, j \in V, i \neq j \tag{2}$$

$$x_{ij}^k \leq y_{ik}, \quad \forall i, j \in V, i \neq j, k \in P(i, j) \tag{3}$$

$$\sum_{l \in \rho(i,j)} \xi_{ij}^l = y_{ij}, \quad \forall i, j \in V, i \neq j \tag{4}$$

$$\sum_{i \in V} \sum_{j \in V} D_{ij}(X) / m_{ij} \sum_{l \in \rho(i,j)} \xi_{ij}^l a_{ij}^{nl} \leq \beta_n C_n^{Link}, \quad \forall n \in E \tag{5}$$

$$\sum_{i \in V} \sum_{j \in V} f_{ij}(X) x_{ij}^k \leq \theta C_R^k, \quad \forall k \in V \tag{6}$$

$$y_{ij} = 1, \quad \forall (i, j) \in S_b \tag{7}$$

$$y_{ij} = 0, \quad \forall (i, j) \in S_f \tag{8}$$

$$\xi_{ij}^l = 1, \quad \forall i, j : l \in \rho_{ij}^S \text{ and } y_{ij} = 1 \tag{9}$$

$$\sum_{j \in V} \phi(D_{kj}(X)) \leq C_{TR}^k, \quad \forall k \in V \tag{10}$$

$$y_{ij}, x_{ij}^k, \xi_{ij}^l \in \{0,1\}. \quad \forall i, j \in V, i \neq j, k \in P(i, j), l \in \rho(i, j) \tag{11}$$

For the meaning of the symbols and constraints of the model, please refer to Lin(2021) for details. Considering the space issue, we don't repeat the description here. The improvements of the model proposed in this paper are expressed as follows：

$$\sum_{i \in V} \sum_{j \in V} D_{ij}(X) / m_{ij} \sum_{l \in \rho(i,j)} \xi_{ij}^l a_{ij}^{nl} \leq \tilde{C}_n^{Link}, \quad \forall n \in E \tag{12}$$

$$\sum_{i \in V} \sum_{j \in V} f_{ij}(X) x_{ij}^k \leq \tilde{C}_R^k, \quad \forall k \in V \tag{13}$$

$$\frac{1}{a} \sum_{j \in V} \phi(D_{kj}(X)) \leq \tilde{C}_{TR}^k, \quad \forall k \in V \tag{14}$$

where $\tilde{C}_n^{Link}$ is the fluctuating available transportation capacity of link $n$, $\tilde{C}_R^k$ is the fluctuating aavailable reclassification capacity of yard $k$ and $\tilde{C}_{TR}^k$ is the fluctuating number of available tracks of yard $k$. The constraint (12) is the elastic link capacity constraint. It indicates that the volume of flows on each arc cannot exceed the fluctuating capacity, the constraint (13) means that the volume of flows reclassified at each yard cannot exceed the fluctuating capacity, and constrain (14) express that the number of occupied tracks is less than the number of fluctuating available classification tracks. In China, it usually occupies one track per 200 cars, the value of $a$ is 200 in this paper.

**3.1 Adaptation to the elastic capacity constraint**

To illustrate the elastic property of the available capacity, we consider a fuzzy version to describe it, which through discussing the membership function that designates the satisfaction degree for the volumes of traffic flows. In fuzzy theory, the determination of the membership function is a key problem (Milenkovie (2013)). For the elastic reclassification capacity constraint of yards, the variable of membership function is the volumes of flows at yards, which is expressed as $\sum_{i \in V}\sum_{j \in V} f_{ij}(X) x_{ij}^k$. In addition, in order to better make membership function designates the satisfaction degree for the volumes of traffic flows, we define an available capacity belt to replace $\tilde{C}_R^k$, we assume that the lower bound of reclassification capacity of yard $k$ is $C_{R\_k}^{lower}$ and that the upper bound is $C_{R\_k}^{upper}$, the available reclassification capacity belt can be formed by the interval $[C_{R\_k}^{lower}, C_{R\_k}^{upper}]$. Let $F_k = \sum_{i \in V}\sum_{j \in V} f_{ij}(X) x_{ij}^k$, the membership function of available reclassification capacities of yards is as follows:

$$\zeta(F_k) = \begin{cases} 1 & F_k \leq C_{R\_k}^{lower} \\ \dfrac{C_{R\_k}^{upper} - F_k}{C_{R\_k}^{upper} - C_{R\_k}^{lower}} & C_{R\_k}^{lower} < F_k \leq C_{R\_k}^{upper} \\ 0 & F_k > C_{R\_k}^{upper} \end{cases} \qquad (15)$$

Through defining the membership function (15), the constraint (13) can be replaced by it. The value of the membership function indicates the satisfaction degree of the reclassification capacity constraint. If the volume of flows reclassified at yard $k$ is smaller than $C_{R\_k}^{lower}$, the value of the function is one, that is, the capacity constraint is satisfied. If the volume of flows at yard $k$ is between $[C_{R\_k}^{lower}, C_{R\_k}^{upper}]$, the value of the function is between (0, 1), that is, the capacity constraint is satisfied with a certain probability. To demonstrate this situation, we use the penalty term to describe it, which is shown in the following content. Finally, if $F_k$ is larger than $C_{R\_k}^{upper}$, the value of function is 0, which means that the capacity constraint is not satisfied.

In order to better describe the effect of $\zeta(F_k)$, we introduce the penalty term to illustrate it. The smaller the value of the membership function, the larger the penalty. For the situation that the value of membership function is between $[C_{R\_k}^{lower}, C_{R\_k}^{upper}]$, the form of penalty term can be described as $\alpha \sum (1 - \zeta(F_k))$, where $\alpha$ is a positive coefficient. In addition, another penalty term is also needed to be added to avoid the volume of shipments exceeds the upper bound of the yard capacity, and this penalty term is $\beta \sum \max(0, F_k - C_{R\_k}^{upper})$. Based on the above analysis, the penalized fuzzy reclassification capacity constraint can be stated as follows:

$$G(X) = \alpha \sum_{k \in V} (1 - \zeta(F_k)) + \beta \sum_{k \in V} \max(0, F_k - C_{R\_k}^{upper}) \qquad (16)$$

In the same way, we assume that the lower bounds of available tracks of yard $k$ and transportation capacity of link $n$ are $C_{TR\_k}^{lower}$ and $C_n^{lower}$, the upper bounds are $C_{TR\_k}^{upper}$ and $C_n^{upper}$ respectively, the available capacity belt of them can be formed by the interval $[C_{TR\_k}^{lower}, C_{TR\_k}^{upper}]$, $[C_n^{lower}, C_n^{upper}]$. Let $T_k^{track} = \dfrac{1}{a}\sum_{j \in V} \phi(D_{kj}(X))$ and

$$R_n^{link} = \sum_{i \in V} \sum_{j \in V} D_{ij} / m_{ij} \sum_{l \in \rho(i,j)} \eta_{ij}^l a_{ij}^{nl}$$, the membership function of available transportation capacities of links and the available tracks of yards are as follows:

$$\vartheta(T_k^{track}) = \begin{cases} 1 & T_k^{track} \leq C_{TR\_k}^{lower} \\ \dfrac{C_{TR\_k}^{upper} - T_k^{track}}{C_{TR\_k}^{upper} - C_{TR\_k}^{lower}} & C_{TR\_k}^{lower} < T_k^{track} \leq C_{TR\_k}^{upper} \\ 0 & T_k^{track} > C_{TR\_k}^{upper} \end{cases} \quad (17)$$

$$\psi(R_n^{link}) = \begin{cases} 1 & R_n^{link} \leq C_n^{lower} \\ \dfrac{C_n^{upper} - R_n^{link}}{C_n^{upper} - C_n^{lower}} & C_n^{lower} < R_n^{link} \leq C_n^{upper} \\ 0 & R_n^{link} > C_n^{upper} \end{cases} \quad (18)$$

The penalized fuzzy transportation capacity constraint of available tracks of yards can be stated as follows: $H(X) = \alpha \sum_{k \in V} (1 - \vartheta(T_k^{track})) + \beta \sum_{k \in V} \max(0, T_k^{track} - C_{TR\_k}^{upper})$, and that of transportation capacity of links is:

$$M(X) = \alpha \sum_{k \in V} (1 - \psi(R_n^{link})) + \beta \sum_{k \in V} \max(0, R_n^{link} - C_n^{upper})$$

It should be noted that the values of penalty factors of the above membership functions are the same, both are $\alpha$ and $\beta$, they take the value of 1500 (Lin (2020)). $G(X)$, $H(X)$ and $M(X)$ are penalized constraints, which are common in heuristic algorithms, especially simulated annealing algorithm.

## 4 Conclusions

This paper studies the integral optimization of train formation plan and traffic routing (TFP&TR), the fuzzy theory is introduced to designate the satisfaction degree for the volumes of traffic flows by defining the membership function. A non-linear integer programming model, is constructed in this paper. In the model, we comprehensively consider various operational requirements and a set of elastic capacity constraints, including link capacity, yard reclassification capacity and the maximal number of blocks a yard can be formed, while trying to minimize the total costs of accumulation, reclassification, and transportation.

Our future research direction includes (1) designing efficient heuristic solution approach to solve the mathematical model; and (2) testing the proposed solution approach and applying the approach to large-scale real-life problem instances.